
\input amstex.tex
\documentstyle{amsppt}

\def\DJ{\leavevmode\setbox0=\hbox{D}\kern0pt\rlap
 {\kern.04em\raise.188\ht0\hbox{-}}D}
\def\dj{\leavevmode
 \setbox0=\hbox{d}\kern0pt\rlap{\kern.215em\raise.46\ht0\hbox{-}}d}

\baselineskip=13pt
\def\hf{{\textstyle{1\over2}}}

\def\d{{\,\roman d}}
\def\e{\varepsilon}

\def\G{\Gamma}
\def\k{\kappa}
\def\s{\sigma}

\def\={\;=\;}

\def\zt{\zeta(\hf+it)}

\def\R{\Re{\roman e}\,}  
\def\z{\zeta}

\def\hf{{\textstyle{1\over2}}}

\font\tenmsb=msbm10
\font\sevenmsb=msbm7
\font\fivemsb=msbm5
\newfam\msbfam
\textfont\msbfam=\tenmsb
\scriptfont\msbfam=\sevenmsb
\scriptscriptfont\msbfam=\fivemsb
\def\Bbb#1{{\fam\msbfam #1}}

\def \NN {\Bbb N}

\def \RR {\Bbb R}

\font\ff=cmr8

\baselineskip=13pt

\font\teneufm=eufm10
\font\seveneufm=eufm7
\font\fiveeufm=eufm5
\newfam\eufmfam
\textfont\eufmfam=\teneufm
\scriptfont\eufmfam=\seveneufm
\scriptscriptfont\eufmfam=\fiveeufm
\def\mathfrak#1{{\fam\eufmfam\relax#1}}

\font\tenmsb=msbm10
\font\sevenmsb=msbm7
\font\fivemsb=msbm5
\newfam\msbfam
     \textfont\msbfam=\tenmsb
      \scriptfont\msbfam=\sevenmsb
      \scriptscriptfont\msbfam=\fivemsb
\def\Bbb#1{{\fam\msbfam #1}}

\def \NN {\Bbb N}

\def \RR {\Bbb R}

  \def\rightheadline{{\hfil{\ff
  Mean value result for the fourth moment of $|\zt|$}\hfil\tenrm\folio}}

  \def\leftheadline{{\tenrm\folio\hfil{\ff
   Aleksandar Ivi\'c }\hfil}}
  \def\emptyheadline{\hfil}
  \headline{\ifnum\pageno=1 \emptyheadline\else
  \ifodd\pageno \rightheadline \else \leftheadline\fi\fi}

\font\ff=cmr8
\font\teneufm=eufm10
\font\seveneufm=eufm7
\font\fiveeufm=eufm5
\newfam\eufmfam
\textfont\eufmfam=\teneufm
\scriptfont\eufmfam=\seveneufm
\scriptscriptfont\eufmfam=\fiveeufm
\def\mathfrak#1{{\fam\eufmfam\relax#1}}

\font\tenmsb=msbm10
\font\sevenmsb=msbm7
\font\fivemsb=msbm5
\newfam\msbfam
\textfont\msbfam=\tenmsb
\scriptfont\msbfam=\sevenmsb
\scriptscriptfont\msbfam=\fivemsb
\def\Bbb#1{{\fam\msbfam #1}}

\def \NN {\Bbb N}

\def \RR {\Bbb R}

 \def\e{\varepsilon}
 \def\d{\,{\roman d}}
\topmatter
\title
A MEAN VALUE RESULT INVOLVING THE FOURTH MOMENT OF $|\zt|$
\endtitle
\author
Aleksandar Ivi\'c
\endauthor
\dedicatory
To Prof. Imre K\'atai on the occasion of his 65th birthday
\enddedicatory
\address
Katedra Matematike RGF-a Universiteta u Beogradu, \DJ u\v sina 7,
11000 Beograd, Serbia (Yugoslavia)
\endaddress
\subjclass 11 M 06
\endsubjclass
\email {\tt  aivic\@rgf.bg.ac.yu, eivica\@ubbg.etf.bg.ac.yu}
\endemail
\dedicatory  
To Prof. Imre K\'atai on the occasion of his 65th birthday
\enddedicatory
\abstract
If $(k,\ell)$ is an exponent pair such that $k + \ell < 1$, then we have
$$
\int_0^T|\zt|^4|\z(\s+it)|^2\d t \ll_\e T^{1+\e}
\quad\left(\s > \min\left({5\over6},
\max\bigl(\ell-k,\,{5k + \ell\over 4k+1}\bigr)\right)\right),
$$
while if $(k,\ell)$ is an exponent pair such that $3k + \ell < 1$, 
then we have
$$
\int_0^T|\zt|^4|\z(\s+it)|^4\d t \ll_\e T^{1+\e}\qquad\left( 
\s > {11k+\ell+1\over8k+2}\right).
$$

\endabstract
\endtopmatter
\head
1. Introduction
\endhead

Let as usual $\z(s) = \sum_{n=1}^\infty n^{-s}\;(\s > 1)$ denote
the Riemann zeta-function, where $s = \s+it$ is a complex
variable. Mean values of $\z(s)$ in the so-called ``critical strip"
$\hf \le \s \le 1$ represent a central topic in the theory of the 
zeta-function (see [7] and [8]   for an extensive account). No bound of
the form
$$
\int_0^T|\zt|^{2m}\d t \;\ll_{\e,m}\; T^{1+\e} \qquad(m \in \NN)\leqno(1.1)
$$
is known to hold when $m\ge3$, while in the cases $m = 1,2$ precise
asymptotic formulas for the integrals in question are known (see op. cit.).
Here we shall prove two hybrid bounds involving the mean value of
$|\zt|^4$ multiplied by $|\z(\s+it)|^{2j}\;(j=1,2;\;\hf<\s<1)$. 
The results are

\medskip
THEOREM 1. 
{\it If $(k,\ell)$ is an exponent pair such that $k + \ell < 1$, 
then we have}
$$
\int_0^T|\zt|^4|\z(\s+it)|^2\d t \ll_\e T^{1+\e}
\quad\left(\s > \min\left({5\over6},
\max\bigl(\ell-k,\,{5k + \ell\over 4k+1}\bigr)\right)\right),
\leqno(1.2)
$$
{\it and in particular} (1.2) {\it holds for} $\s \ge 5/6 = 0.8333\ldots\,$.

\medskip
THEOREM 2. 
{\it If $(k,\ell)$ is an exponent pair such that $3k + \ell < 1$, 
then we have}
$$
\int_0^T|\zt|^4|\z(\s+it)|^4\d t \ll_\e T^{1+\e}
\qquad\left(\s > \max \left({\ell-k+1\over2},\,
{11k+\ell+1\over8k+2}\right)\right),\leqno(1.3)
$$
{\it and in particular} (1.3) {\it holds for $\s \ge 1953/1984 =
0.984375$.}

\medskip The merit of these results is that (1.2) and (1.3) hold for
values of $\s$ less than one; of course one expects the bounds to
hold for $\s \ge \hf$, in which case we would obtain the (yet unproved)
sixth and eighth moment of $|\zt|$ (namely (1.1) with $m = 3$ and $m = 4$,
respectively).

\head
2. Proof of Theorem 1
\endhead
In the proof of both (1.2) and (1.3) it is sufficient to consider the
integral over $[T, 2T]$, then to replace $T$ by $T2^{-j}$ $(j = 1,2,\ldots)$
and sum all the resulting estimates. Also, it is sufficient to suppose
that $\s \le 1$, since one has (see e.g., [7])
$$
\z(\s+it) \;\ll\;\log|t|\qquad(\s \ge 1).
$$
To prove the bound on $\s$ in (1.2) involving $k,\ell$, we shall
use the simple approximate
functional equation for $\z(s)$ (see [7, Theorem 1.8]), which gives
$$
\z(s) = \sum_{n\le T}n^{-s} + O(1)\qquad(s = \s+it,\,T \le t \le 2T).
\leqno(2.1)
$$
The essential tool in our considerations is the following theorem for
the fourth moment of $|\zt|$, weighted by a Dirichlet polynomial, due
to N. Watt [9]. This is built on the works of J.-M. Deshouillers 
and H. Iwaniec [1], [2], involving the use of Kloosterman sums,
but it contains the following sharper result: Let $a_1, a_2,
\ldots$ be complex numbers. Then, for $\e > 0, M\ge 1$ and $T\ge 1$,
$$
\int_0^T|\sum_{m\le M}a_m m^{it}|^2|\zt|^4\d t \ll_\e
T^{1+\e}M(1 + M^2T^{-1/2})\max_{m\le M}|a_m|^2.\leqno(2.2)
$$
Here and later $\e$ denotes arbitrarily small, positive constants, not
necessarily the same ones at each occurrence. We write (2.1) as
$$\eqalign{
\z(s) &= \sum_{m\le Y}m^{-s} + \sum_{Y<n\le T}n^{-s} + O(1)\cr&
= \sum\nolimits_1 + \sum\nolimits_2 + \,O(1),\cr}\leqno(2.3)
$$
say, where $1 \ll Y \le T$.
The sum $\sum_1$ is split into $O(\log T)$ subsums with
$N < m \le N' \le 2N\le Y$. In (2.1) we take 
$a_m = m^{-\s}$ for $N < m \le N'$,
$a_m = 0$ otherwise. Then in view of $N \ll Y, \hf \le \s \le 1$ it
follows that
$$\eqalign{
\int_0^T|\sum\nolimits_1|^2|\zt|^4\d t &\ll_\e T^{1+\e}\max_{N\ll Y}
N^{1-2\s}(1 + N^2T^{-1/2}) \cr&
\ll_\e T^{1+\e}(1 + Y^{3-2\s}T^{-1/2}) \ll_\e T^{1+\e}\cr}\leqno(2.4)
$$
for
$$
Y \;=\; T^{1\over6-4\s}.\leqno(2.5)
$$
To estimate $\sum_2$ in (2.3) we use the theory of (one-dimensional)
exponent pairs (see e.g., [3], [5] and [7]). We split $\sum_2$ into 
$O(\log T)$ subsums with $N < m \le N' \le 2N$, $Y \le N \le T$ and
$\s\ge\hf$. Removing the (monotonically decreasing) factor $n^{-\s}$
by partial summation from each subsum, it remains to estimate
$$
S(N,t) := \sum_{N<n\le N'\le2N}n^{it}\qquad(Y \le N \le T,\,
T\le t \le 2T).
$$
If $(k,\ell)$ is an exponent pair, then since $n^{it} = {\roman e}^{iF(n,t)}$
with ${\partial^r F(n,t)\over\partial n^r} \asymp_r TN^{-r}$, it follows
that 
$$
S(N,t) \;\ll\;\left({T\over N}\right)^kN^\ell = T^kN^{\ell-k},
$$
and consequently
$$
\sum\nolimits_2 \;\ll\;T^kN^{\ell-k-\s}\log T \ll 
T^{k+{\ell-k-\s\over6-4\s}}\log T
$$
if $\s \ge \ell-k$, which is  our assumption.
Hence $\sum_2 \ll \log T$ for $k+{\ell-k-\s\over6-4\s}\le0$,
i.e.
$$
\s \;\ge\;{5k+\ell\over4k+1},
$$
giving
$$
\int_0^T|\sum\nolimits_2|^2|\zt|^4\d t 
\; \ll_\e \;T^{1+\e}\quad\left(\s \ge 
\max\Bigl(\ell-k,\,{5k + \ell\over4k+1}\Bigr)\right).\leqno(2.6)
$$
Combining (2.4) and (2.6) we obtain the second bound in 
(1.2); for $k+\ell < 1$
we have ${5k + \ell\over4k+1} < 1$. To obtain a specific result
we choose M.N. Huxley's exponent pair (see [6]) $(\k,\lambda) =
({32\over205}+\e,\,\hf + {32\over205}+\e)$, which supersedes his
exponent pair (see [4], [5]) $(\k,\lambda) =
({89\over570}+\e,\,{374\over570}+\e)$. This exponent pair
is one of the many obtained by the Bombieri--Iwaniec method.
With this pair we find that
$$
{5k + \ell\over4k+1} \;=\;{589\over666} \;=\; 0.884384384\ldots.
\leqno(2.7)
$$
As is often the case when
one applies the theory of exponent pairs, the above exponent pair
is not optimal, and small improvements may be obtained by more
laborious calculations. Note that the algorithm of Graham-Kolesnik 
[3, Chapter 5] cannot be used when the exponent pairs are formed
by the use of (variants) of the Bombieri--Iwaniec method, and not only
by the classical $A-, B$-process and convexity, so the optimal value
is hard to compute. However, in the above case the exponent pair is
in a certain sense optimal. Namely if $\ell = k + \hf$, then one has
(see [3, Theorem 4.1])
$$
\mu(\hf) \le k, \quad \mu(\s) := \limsup_{t\to\infty}\,{\log|\z(\s+it)|
\over\log t}\quad(\s\in\RR).
$$
But if $\ell = k + \hf$, then
$$
{5k + \ell\over4k+1} = {\hf + 6k\over1+4k},
$$
which is an upper bound for
$$
{\hf + 6\mu(\hf)\over1+4\mu(\hf)}.\leqno(2.8)
$$
If we use the bound ([7, eq. (8.14)])
$$
\sum_{N<n\le2N}n^{-\s-it} \ll_\e N^{\s_0-\s}T^{\e-\s_0}
+ N^{\s_0-\s}\int_0^T|\z(\s_0 + it+iv)|\,{\d v\over v +1},
$$
then we obtain
$$
\sum_{N<n\le2N}n^{-\s-it} \ll_\e 1 + N^{{1\over2}-\s}T^{\mu({1\over2})+\e}
\quad(T\le t \le 2T,\,\s > \hf).
$$
But if $N \ge Y = T^{1/(6-4\s)}$, then the above bound gives
$$
\sum_{N<n\le2N}n^{-\s-it} \ll_\e T^\e\quad(T\le t \le 2T,\,N \ge Y,\,\s > \hf)
$$
for 
$$
\s \ge {\hf + 6\mu(\hf)\over1+4\mu(\hf)},\leqno(2.9)
$$
which is (2.8). Huxley's work [6] brings forth precisely
the new bound (hitherto the sharpest one of its kind)
$\mu(\hf) \le 32/205$, corresponding to the exponent pair
$(k,\ell)$ with $k = 32/205+\e$, $\ell = k + \hf$, so that
in this context the value given by (2.7) is the optimal one that
can be obtained at present from exponent pairs satisfying the
condition $\ell = k + \hf$.

\medskip
To complete the proof of (1.2) we use the well-known Mellin inversion
integral
$$
{\roman e}^{-x} = {1\over2\pi i}\int_{(c)}x^{-w}\G(w)\d w
\quad(c>0,\,x>0),
\leqno(2.10)
$$
where $\int_{(c)}$ denotes integration over the line $\R w = c$. In
(2.10) we set $x = n/Y\;(1 \ll Y \ll T^C)$, multiply by $n^{-s}\;(\hf <
\s < 1)$  and sum over $n$. This gives
$$
\sum_{n=1}^\infty {\roman e}^{-n/Y}n^{-s} =  
{1\over2\pi i}\int_{(2)}Y^w\z(s+w)\G(w)\d w
\quad(s = \s+it,\,T \le t \le 2T).
\leqno(2.11)
$$
In (2.11) we shift the line of integration to $\R w = \hf -\s$ and
apply the residue theorem. The pole at $w = 1 -s$ contributes a residue
which is, by Stirling's formula for $\G(s)$, $\ll 1$. The pole at
$w = 0$ yields $\z(s)$, and we obtain from (2.11)
$$
\z(s) \ll 1 + \left|\sum_{n\le Y\log^2T}{\roman e}^{-n/Y}n^{-\s-it}\right|
+ Y^{{1\over2}-\s}\int_{-\log^2T}^{\log^2T}|\z(\hf + it + iv)|\d v.
$$
Therefore 
$$
\int_0^T|\zt|^4|\z(\s+it)|^2\d t \ll T\log^4T + I_1(T) + I_2(T),\leqno(2.12)
$$
say, where
$$\eqalign{
I_1(T) :&= \int_0^T|\zt|^4
\Bigl|\sum_{n\le Y\log^2T}{\roman e}^{-n/Y}n^{-\s-it}\Bigr|^2\d t,\cr
I_2(T) :&= Y^{1-2\s}\int_0^T|\zt|^4
\left(\int_{-\log^2T}^{\log^2T}|\z(\hf + it + iv)|\d v\right)^2\d t.\cr}
\leqno(2.13)
$$
Similarly to (2.4) we obtain
$$
I_1(T) \;\ll_\e\; T^{1+\e}(1 + Y^{3-2\s}T^{-1/2}).\leqno(2.14)
$$
To $I_2(T)$ we apply H\"older's inequality for integrals and the sharpest
bound for the sixth moment of $|\zt|$ (see [7, Chapter 8]),
namely $\int_0^T|\zt|^6\d t \ll T^{5/4}\log^CT$, to deduce that
$$
I_2(T)  \;\ll_\e\; T^{{5\over4}+\e}Y^{1-2\s}.\leqno(2.15)
$$
Now we choose
$$
Y \;=\;T^{3/8}.
$$
Then from (2.12)--(2.15) it follows that
$$
\int_0^T|\zt|^4|\z(\s+it)|^2\d t 
\ll_\e T^{1+\e} + T^{{1\over2}+{9-6\s\over8}+\e} \ll_\e T^{1+\e}
$$
for $\s \ge 5/6$, which yields the first bound in (1.2) and completes
the proof of Theorem 1.

\head
3. Proof of Theorem 2
\endhead
For the proof of Theorem 2 we shall use the approximate functional
equation (see [7, Theorem 4.2])
$$
\z^2(s) = \sum_{n\le x}d(n)n^{-s} + \chi^2(s)\sum_{n\le y}d(n)n^{s-1}
+ O(x^{{1\over2}-\s}\log t),\leqno(3.1)
$$
where $d(n)$ is the number of divisors of $n$, $0<\s<1;\, x,y,t > C > 0$
and $4\pi^2xy = t$. Here
$$
\chi(s) = {\z(s)\over\z(1-s)} = 2^s\pi^{s-1}\G(1-s)\sin(\hf\pi s)
\asymp t^{{1\over2}-\s}\quad(t\ge t_0 > 0)
$$
is the expression appearing in the functional equation for $\z(s)$. In
(3.1) we suppose that $T \le t \le 2T, T \le x \le 2T$. Then we obtain
$$\eqalign{
|\z(\s+it)|^4 & \ll \log^2T +\cr&
+ {\bigl|\sum_{n\le x}d(n)n^{-s}\bigr|}^2 + 
T^{1-2\s}{\bigl|\sum_{n\le {t^2\over4\pi^2x}}d(n)n^{s-1}\bigr|}^2 .\cr}
$$
Both sums on the right-hand side are split into $O(\log T)$ subsums
with $N < n \le N'\le 2N,\,N\ll T$. Setting 
$S(u) := \sum_{N<n\le u}d(n)n^{-it}$ we have, by partial summation,
$$\eqalign{
\sum_{N<n\le N'}d(n)n^{-\s-it}  &= S(N')(N')^{-\s} + \s\int_N^{N'}
S(u)u^{-\s-1}\d u,\cr
\sum_{N<n\le N'}d(n)n^{\s-1-it}  &= S(N')(N')^{\s-1} + \s\int_N^{N'}
S(u)u^{\s-2}\d u.\cr}
$$
This gives
$$\eqalign{&
{\bigl|\sum_{n\le x}d(n)n^{-s}\bigr|}^2  + 
T^{1-2\s}{\bigl|\sum_{n\le {t^2\over4\pi^2x}}d(n)n^{s-1}\bigr|}^2\cr&
\ll \log T\max_N\left(N^{-2\s}\max_{N\le u\le N'}|S(u)|^2\right)
\left(1 + \left(\frac{T}{N}\right)^{1-2\s}\right)\cr&
\ll \log T\max_N N^{-2\s}\max_{N\le u\le N'}|S(u)|^2,\cr}
$$
since $N\ll T, \s\ge\hf$. In the case when $N \le Y$ (see (2.5)) we have,
by (2.2),
$$\eqalign{&
\int_T^{2T}N^{-2\s}|S(u)|^2|\zt|^4\d t \cr&\ll_\e T^{1+\e}N^{1-2\s}
(1 + N^2T^{-1/2})\max_{N<n\le N'}d^2(n)\cr&
\ll_\e T^{1+\e}(1 + Y^{3-2\s}T^{-1/2}) \ll_\e T^{1+\e},\cr}
$$
since $d(n) \ll_\e n^\e$. 

\smallskip
In the case when $Y < N \ll T$ we shall estimate
$$
{\bar S(u)} = \sum_{N<n\le u}d(n)n^{it} = \sum_{n\le u}d(n)n^{it} 
- \sum_{n\le N}d(n)n^{it} 
$$
by estimating
$$
\sum(u,t) \;:=\; \sum_{n\le u}d(n)n^{it}\qquad(N < u \le N'\le 2N).
$$
By applying the familiar hyperbola method we have
$$\eqalign{&
\sum(u,t) = \sum_{mn\le u}(mn)^{it} \cr&
= 2\sum_{m\le\sqrt{u}}m^{it}\,\sum_{n\le u/m}n^{it} - \left(
\sum_{m\le\sqrt{u}}m^{it}\right)^2\cr&
= 2S_1(u,t) - S_2^2(u,t),\cr}
$$
say. To estimate $S_1(u,t)$, we split the inner sum over $n$ into 
$O(\log T)$ subsums
$$
S_3(u_1,t) := \sum_{u_1<n\le u_1'\le2u_1}n^{it}\qquad(u_1 \ll u/m).
$$
Then, since $\ell \ge k$ for any exponent pair $(k,\ell)$,
$$
S_3(u,t) \;\ll\; T^k\left(\frac{u}{m}\right)^{\ell-k},
$$
which yields
$$
S_1(u,t) \;\ll\; T^k\log T\cdot N^{\ell-k}\sum_{m\ll\sqrt{N}}m^{k-\ell}
\ll T^kN^{{1\over2}(\ell-k+1)}\log T.
$$
In a similar vein it follows that
$$
S_2(u,t) \ll T^kN^{{1\over2}(\ell-k)}\log T,
$$
and thus for $N \gg Y$ we obtain, for $\s \ge \hf(\ell-k+1)$, 
$$
\eqalign{
N^{-\s}|S(u)| &
\ll N^{-\s}\log^2T(T^kN^{{1\over2}(\ell-k+1)} + T^{2k}N^{\ell-k})\cr&
\ll (T^kT^{\ell-k+1-2\s\over12-8\s} + 
T^{2k}T^{\ell-k-\s\over6-4\s})\log^2T\cr&
\ll \log^2T\cr}
$$
for
$$
\s \;\ge\; \max\left({11k+\ell+1\over8k+2}, {11k+\ell\over8k +1}\right)
= {11k+\ell+1\over8k+2}
$$
if $3k +\ell < 1$, which we supposed.
This proves (1.3). Finally we consider the exponent pair (see [3, p. 39])
$$
(k,\ell) = \left({16\over120Q-32},\,{120Q-16q-63\over120Q-32}\right)
\quad (Q = 2^q,\; q\ge 2).
$$
The optimal value for $q$ is in our case found to be $q=3$, giving
$$
(k,\ell) = \left({16\over928},\,{849\over928}\right)
= {11k+\ell+1\over8k+2} \;=\; {1953\over1984} = 0.984375
\;\left( > \hf(\ell-k+1)\right).\leqno(3.2)
$$
This completes the proof of Theorem 2, and with a more careful choice
of the exponent pair the value (3.2) could be improved a little
(namely by the use of the algorithm of [3, Chapter 5]). 
It is an open problem to find  $\s_0 = \s_0(j) \;(<1)$ such that
$$
\int_0^T|\zt|^4|\z(\s+it)|^{2j}\d t \ll_{j,\e} T^{1+\e}
$$
for $j\in \NN$ satisfying $j\ge3$ and $\s > \s_0$. By using the method
outlined at the end of Section 2, one would obtain the value
$$
\s_0 \le {\hf + 6j\mu(\hf)\over1 + 4j\mu(\hf)}.
$$
But the right-hand side does not exceed unity if and only if
$$
\mu(\hf) \;\le\;{1\over4j},
$$
which is not known to hold unless $j=1$, and this case we already
considered. Thus the approach based on the use of exponent pairs
seems more appropriate already in the case $j=2$.

\vfill
\eject\topskip2cm

\bigskip\bigskip
\Refs
\bigskip
\item{[1]} J.-M. Deshouillers and H. Iwaniec, Power mean-values of the
the Riemann zeta-function, Mathematika {\bf29}(1982), 202-212.

\item{[2]} J.-M. Deshouillers and H. Iwaniec, Kloosterman sums and
Fourier coefficients of cusp forms,  Invent. Math. {\bf70}(1982), 219-288.

\item{[3]} S.W. Graham and G. Kolesnik, Van der Corput's method
of exponential sums, LMS Lecture Notes series {\bf126}, Cambridge
University Press, Cambridge, 1991.

\item{[4]} M.N. Huxley, Exponential sums and the Riemann zeta function
IV, Proc. London Math. Soc. (3){\bf60} (1993), 1-40.

\item{[5]} M.N. Huxley, Area, lattice points and exponential sums,
LMS Monographs (New Series) {\bf13}, Oxford University Press,
Oxford, 1996.

\item{[6]} M.N. Huxley, Integer points, exponential sums and the Riemann
zeta function, in ``Number Theory for the Millenium. Proc. Millenial Conf.
on Number Theory" (Urbana, 2000),
in print.

\item{[7]} A. Ivi\'c, The Riemann zeta-function, John Wiley \&
Sons, New York, 1985.

\item{[8]} A. Ivi\'c, The mean values of the Riemann zeta-function, Tata 
	Institute of Fundamental Research, Lecture Notes {\bf82}, 
    Bombay 1991 (distr. Springer Verlag, Berlin etc.).

\item{[9]} N. Watt, Kloosterman sums and a mean value for Dirichlet
polynomials, J. Number Theory {\bf53}(1995), 179-210.

\bigskip

Aleksandar Ivi\'c

Katedra Matematike RGF-a

Universitet u Beogradu

\DJ u\v sina 7, 11000 Beograd, Serbia

{\tt aivic\@rgf.bg.ac.yu}

\endRefs


\bye